\documentclass[12pt]{article}
\usepackage{amsmath} 
\usepackage{graphicx}
\usepackage[all]{xy}
\usepackage{amssymb}
\input diagxy
\usepackage[usenames]{color}

\newtheorem{prop}{Proposition}[section]

\newtheorem{rem}{Remark}[section]
\newtheorem{ex}{Example}[section]

\newtheorem{lem}{Lemma}[section]
\newtheorem{dfn}{Definition}[section]

\def\mathopdef#1{\expandafter\def\csname#1\endcsname{\mathop{\rm#1}\nolimits}}
\def\mathoplsdef#1{\expandafter\def\csname#1\endcsname{\mathop{\rm#1}}}
\def\mathbfdef#1{\expandafter\def\csname#1\endcsname{{\rm\bf#1}}}
\def\mathrmdef#1{\expandafter\def\csname#1\endcsname{{\rm#1}}}

\newcommand{\pf}{\vspace{\baselineskip}\noindent{\bf Proof. }}
\newcommand{\qed}{\hfill\rule{4pt}{8pt}\par\vspace{\baselineskip}}

\newcommand{\bA}{{\bf A}}

\newcommand{\bG}{{\bf G}}

\newcommand{\ba}{{\bf a}}
\newcommand{\bb}{{\bf b}}
\newcommand{\bc}{{\bf c}}

\newcommand{\BB}{\mathbf{BB}}

\newcommand{\tw}{{\tau}}

\newcommand{\rel}{\mathbf{Rel}}

\newcommand{\tr}{\mathbf{TRel}}

\title{Racks and blocked braids}
\author{
N. Sabadini\\
Dipartimento di Scienza e Alta Tecnologia\\
Universit\`a dell' Insubria\\
via Carloni, 78, Como, Italy\\
{\tt nicoletta.sabadini@uninsubria.it}
\and
R. F. C. Walters\\
Dipartimento di Scienza e Alta Tecnologia\\
Universit\`a dell' Insubria\\
via Carloni, 78, Como, Italy\\
{\tt robert.walters@uninsubria.it}}
\begin{document}
\maketitle
\begin{abstract}
In the paper \cite{MSW13} the present authors together with Davide Maglia introduced the blocked-braid groups $\BB_n$ on $n$ strands, and proved that a blocked torsion has order either $2$ or $4$. We conjectured that the order was actually $4$ but our methods in that paper, which involved introducing for any group $\bG$ a braided monoidal category of \emph{tangled relations}, were inadequate to demonstrate this fact.  
Subsequently Davide Maglia in unpublished work investigated exactly what part of the structure and properties of a group $\bG$ are needed to permit the construction of a braided monoidal category with a tangle algebra and was able to distinguish blocked two-torsions from the identity.

In this paper we present a simplification of his answer, which turns out to be related to the notion of rack. We show that if $\bG$ is a rack then there is a 
braided monoidal category  $\tr_{\bG}$ generalizing that of \cite{MSW13}. Further we introduce a variation of the notion of rack which we call \emph{irack} which yields a tangle algebra in $\tr_{\bG}$. Iracks are in particular racks but have in addition to the operations  abstracting group conjugation also a unary operation abstracting group inverse. Using iracks we obtain new invariants for tangles and blocked braids permitting us to give a proof of Maglia's result that
a blocked double torsion is not the identity.

This work was presented at the Conference in Memory of Aurelio Carboni, Milan, 24-26 June 2013.
\end{abstract}


\section{Racks and Iracks}
We recall the definition of rack (\cite{W,CW59,J82,RF92}).
\begin{dfn}
A \emph{rack} $\bA$ consists of a set $A$ with two binary operations $\lhd$ and $\rhd$ satisfying the equations:
\begin{itemize}
\item[R(1)] $ (a\rhd b)\lhd a=b$,
\item[R(2)] $ a\rhd (b\lhd a)=b $,
\item[R(3)] $a\rhd (b\rhd c)=(a\rhd b)\rhd (a\rhd c)$,
\item[R(4)] $(c\lhd b)\lhd a=(c\lhd a)\lhd (b\lhd a)$.
\end{itemize}
\end{dfn}
The new notion of \emph{irack} which we introduce in this paper is defined as follows:
\begin{dfn}
A \emph{irack} $\bA$ consists of a set $A$ with two unary operations $()^+:A\to A$, $()^-:A\to A$ and 
a binary operation $\rhd :A\times A\to A$ satisfying the equations:
\begin{itemize}
\item[IR(1)] $(a^+)^-=a=(a^-)^+$,
\item[IR(2)] $a\rhd a^-=a^+ $,
\item[IR(3)] $a^-\rhd (a\rhd b)=b $,
\item[IR(4)] $a\rhd (a^-\rhd b)=b $,
\item[IR(5)] $a\rhd(b\rhd c)=(a\rhd b)\rhd (a\rhd c) $,
\item[IR(6)] $ (a\rhd b^-)=(a\rhd b)^-$.
\end{itemize}
We denote the irack $\bA$ as $\bA=(A,\rhd,+,-)$.
\end{dfn}

\begin{rem}
Irack axioms (3), (4), (5) and (6) imply that an irack is a rack if we define $b\lhd a$ to be $a^-\rhd b$. 
\end{rem}

\begin{ex}
Given a group $\bG$ if we define we define $g\rhd h$ to be $ghg^{-1}$ and $h\lhd g$ to be $g^{-1}hg$ then $\bG$ is a rack.  If we define further $g^-$ and $g^+$ both to be the inverse of $g$ then $\bG$ is an irack. 
\end{ex}

\begin{ex}
Given $A=\{1,a,b,c,d,e,f\}$ the following define an irack structure on $A$:

\centerline{
\begin{tabular}{|l|l|l|l|l|l|l|l|}
\hline
       & 1 & a & b & c & d & e & f \\
\hline
$()^+$ & 1 & b & a & f & c & d & e  \\
\hline
$()^-$ & 1 & b & a & d & e & f & c  \\
\hline
\end{tabular}
}

\bigskip
\centerline{
\begin{tabular}{|l|l|l|l|l|l|l|l|}
\hline
$\rhd$ & 1 & a & b & c & d & e & f \\\hline
1      & 1 & a & b & c & d & e & f \\\hline
a      & 1 & a & b & c & d & e & f \\\hline
b      & 1 & a & b & c & d & e & f \\\hline
c      & 1 & b & a & e & f & c & d \\\hline
d      & 1 & b & a & e & f & c & d \\\hline
e      & 1 & b & a & e & f & c & d \\\hline
f      & 1 & b & a & e & f & c & d \\\hline
\end{tabular}
}
\end{ex}

\noindent 
The following lemmas hold for elements of an irack $\bA$.

\begin{lem}
$$a\rhd a=a^{++}.$$
\end{lem}

\pf
\begin{align}
a\rhd a&=(a\rhd a)^{--++}\tag{by IR(1)}\\
&=(a\rhd a^{-})^{-++}\tag{by IR(6)}\\
&=a^{+-++}\tag{by IR(2)}\\
&=a^{++} \tag{by IR(1)}
\end{align}
\qed

\begin{lem}
$$a^-\rhd a=a^{++}.$$
\end{lem}

\pf
\begin{align}
a^-\rhd a&=(a^-\rhd a)^{--++}\tag{by IR(1)}\\
&=(a^-\rhd a^{--})^{++}\tag{by IR(6)}\\
&=a^{-+++}\tag{by IR(2)}\\
&=a^{++} \tag{by IR(1)}
\end{align}
\qed

\begin{lem}
$$a^{++++}=a.$$
\end{lem}

\pf
\begin{align}
(a^{++})^{++}&=(a\rhd a)^{++}\tag{by lemma 1.1}\\
&=(a\rhd a)\rhd (a\rhd a)\tag{by lemma 1.1}\\
&=a\rhd (a\rhd a)\tag{by IR(5)}\\
&=a\rhd a^{++} \tag{by lemma 1.1}\\
&=a\rhd (a^-\rhd a) \tag{by lemma 1.2}\\
&=a \tag{by IR(4)}
\end{align}
\qed

\begin{lem}
$$a^{--}\rhd a=a^{++}.$$
\end{lem}

\pf
\begin{align}
(a^{--}\rhd a)&=(a^{--}\rhd a)^{--++}\tag{by IR(1)}\\
&=(a^{--}\rhd a^{--})^{++}\tag{by IR(6)}\\
&=a^{--++++} \tag{by lemma 1.1}\\
&=a^{++} \tag{by IR(1)}
\end{align}
\qed

\begin{lem}
$$(a\rhd b)^+=a\rhd b^+.$$
\end{lem}

\pf
\begin{align}
(a\rhd b)^+&=(a\rhd b)^{-++}\tag{by IR(1)}\\
&=(a\rhd b^{-})^{++}\tag{by IR(6)}\\
&=(a\rhd b^{-})^{--}\tag{by lemma 1.3}\\
&=(a\rhd b^{---}) \tag{by IR(6)}\\
&=(a\rhd b^{-++}) \tag{by lemma 1.3}\\
&=(a\rhd b^{+}) \tag{by IR(1)}
\end{align}
\qed

\begin{lem}
$$a\rhd a^+=a^-.$$
\end{lem}

\pf
\begin{align}
a\rhd a^+&=(a\rhd a)^{+}\tag{by lemma 1.5}\\
&=a^{+++}\tag{by lemma 1.1}\\
&=a^{--+}\tag{by lemma 1.3}\\
&=a^{-} \tag{by IR(1)}
\end{align}
\qed

\begin{lem}
$$a^+\rhd (a\rhd b)=b.$$
\end{lem}

\pf
\begin{align}
a^+\rhd (a\rhd b)&=a^+\rhd ((a^{+})^-\rhd b)\tag{by IR(1)}\\
&=b\tag{by IR(4)}
\end{align}
\qed

\begin{lem}
$$a\rhd (a^+\rhd b)=b.$$
\end{lem}

\pf
\begin{align}
a\rhd (a^+\rhd b)&=(a^+)^-\rhd (a^{+}\rhd b)\tag{by IR(1)}\\
&=b\tag{by IR(3)}
\end{align}
\qed

\begin{lem}
The  function $x\mapsto a\rhd x$ from $A$ to $A$ has inverse $x\mapsto a^-\rhd x$.
\end{lem}

\pf
This follows directly from properties IR(3) and IR(4).
\qed

\begin{lem}
$$a^+\rhd b=a^-\rhd b.$$
\end{lem}

\pf
Lemmas 1.7 and 1.8 imply that the function $x\mapsto a^+\rhd x$ is also inverse to $x\mapsto a\rhd x$
and hence is equal to $x\mapsto a^-\rhd x$. 
\qed

The lemmas above clearly imply the following proposition:
\begin{prop}
If $\bA=(A,\rhd,+,-)$ is an irack then so also is $\bA=(A,\rhd,-,+)$.
\end{prop}

\section{Racks of $n$ tuples}

\begin{dfn}
If $\bA = (A,\rhd,\lhd)$ is a rack, $\ba=(a_1,a_2,\cdots ,a_n)$ is an $n$-tuple of elements of $A$ (for some $n$, $n=0,1,2,\cdots$) and $b$  is an element of $A$ then we define $\ba \rhd b$ as $\ba \rhd b =b$ if $n=0$ and otherwise
$$\ba \rhd b=a_1\rhd (a_2\rhd (a_3\rhd (\cdots (a_n\rhd b)\cdots ))).$$ 
Similarly we define $b \lhd \ba$ as $b\lhd \ba =b$ if $n=0$ and otherwise
$$b\lhd \ba =      ( ((\cdots (b\lhd a_1)\cdots )\lhd a_{n-2})\lhd a_{n-1})\lhd a_n .$$ 

If also 
$\bb$ is an m-tuple of elements of $\bA$ (for some $m$, $m=0,1,2,\cdots$) we define $\ba\rhd\bb$ as $\ba\rhd\bb=()$ if $m=0$ else as
$$\ba\rhd (b_1,b_2,\cdots ,b_m)=(\ba\rhd b_1,\ba\rhd b_2,\cdots ,\ba\rhd b_m).$$ 
Similarly, we define $\bb\rhd\ba$ as $\bb\rhd\ba=()$ if $m=0$ else as
$$(b_1,b_2,\cdots ,b_m)\lhd \ba=(b_1\lhd \ba,b_2\lhd \ba,\cdots ,b_m\lhd \ba).$$ 

\end{dfn}

\begin{lem}
If $\ba$ is an n-tuple of elements of a rack and $b$ and $c$ are elements of the rack then
$$\ba\rhd (b\rhd c)=(\ba\rhd b)\rhd(\ba\rhd  c).$$
\end{lem}

\pf
\begin{align}
\ba\rhd (b\rhd c) &= a_1\rhd ((a_2,\cdots ,a_n)\rhd (b\rhd c)) \tag{by definition}\\
&=a_1\rhd (((a_2,\cdots ,a_n)\rhd b)\rhd((a_2,\cdots ,a_n)\rhd  c))\tag{by inductive hypothesis)}\\
&=(a_1\rhd ((a_2,\cdots ,a_n)\rhd b))\rhd(a_1\rhd ((a_2,\cdots ,a_n)\rhd  c))\tag{by R(3)}\\
&=(\ba\rhd b)\rhd(\ba\rhd  c) \tag{by definition}
\end{align}
\qed

\begin{lem}
If $\ba$, $\bb$ and $\bc$ are tuples of elements of a rack then
$$\ba\rhd (\bb\rhd \bc)=(\ba\rhd \bb)\rhd(\ba\rhd  \bc).$$
\end{lem}

\pf
It clearly suffices to consider the case in which $\bc=c$, a single element. 
\begin{align}
\ba\rhd ((b_1,\cdots ,b_m)\rhd c) &= \ba\rhd (b_1\rhd ((b_2,\cdots ,b_m)\rhd c)) \tag{by definition}\\
&=  (\ba\rhd b_1)\rhd (\ba\rhd ((b_2,\cdots ,b_m)\rhd c))\tag{by R(3)}\\
&=(\ba\rhd b_1)\rhd ((\ba\rhd (b_2,\cdots ,b_m))\rhd(\ba\rhd c))\tag{by inductive hypothesis}\\
&=(\ba\rhd b_1)\rhd ((\ba\rhd b_2,\cdots ,\ba\rhd b_m)\rhd(\ba\rhd c)) \tag{by definition}\\
&= (\ba\rhd b_1,\cdots ,\ba\rhd b_m)\rhd(\ba\rhd c) \tag{by definition}\\
&= (\ba\rhd \bb)\rhd(\ba\rhd c) \tag{by definition}
\end{align}

\qed

\begin{lem}
If $\ba$ and $\bb$ are tuples of elements of a rack then
$$\ba\rhd (\bb\lhd \ba)=\bb.$$
\end{lem}

\pf
It clearly suffices to consider the case in which $\bb=b$, a single element. 
\begin{align}
\ba\rhd (b\lhd \ba) &= a_1\rhd (\cdots (a_{n-1}\rhd ( a_n\rhd
                       [((\cdots (b\lhd a_1)\cdots )\lhd a_{n-1})\lhd a_n]
                     ))) \tag{by definition}\\
&=  a_1\rhd (\cdots (a_{n-1}\rhd 
                      [(\cdots (b\lhd a_1)\cdots )\lhd a_{n-1}]
                    )) \tag{by R(2)}\\
&=b\tag{by inductive hypothesis}
\end{align}

\qed





We have now verified that tuples of a rack satisfy the equations R(2) and R(3) of racks; R(1) and R(4) follow similarly and hence we have proved the following result:
\begin{prop}
If $\bA=(A,\rhd,\lhd)$ is a rack then the set of tuples of elements of $\bA$ form a rack, with the definition of the operations as given above.
\end{prop}

\begin{rem}
Notice that we have shown that the set of tuples of an irack has a rack structure. 
We may try to add an irack structure defining $\ba^-$ and $\ba^+$  as $()$ if $n=0$ else as
$$\ba^- = (a_n^-,a_{n-1}^-,\cdots ,a_1^-),\,\,\, \ba^+ = (a_n^+,a_{n-1}^+,\cdots ,a_1^+). $$
With these candidates for the operations $()^-$ and $()^+$ the set of tuples of an irack satisfy all 
the axioms of an irack except axiom IR(2). In the irack of example 1.2 $(a,c)\rhd (a,c)^-=(f,a)$ whereas $(a,c)^+=(f,b) $.
\end{rem}

\section{Tangled relations}
In this section we define a category $\tr_{\bA}$ of tangled relations over an irack $\bA$, extending the case
of tangled relations over a group introduced in \cite{RSW12}. Notice that in defining the braided monoidal structure 
on $\tr_{\bA}$ we need only that $\bA$ is a rack.

\bigskip
\noindent{\bf Note well:} This category will \emph{not} provide invariants for tangled circuits in general, which are arrows in a free braided strict monoidal category on a monoidal graph, whose objects have with \emph{commutative Frobenius algebra structures}. We will show that the category $\tr_{\bA}$ is generated by a tangle algebra object, and hence yields invariants for arrows in a free braided strict monoidal category on a monoidal graph, whose objects have with \emph{tangle algebra structures}.  The blocked braid groups were defined in this context.

\begin{dfn}
Consider an irack $\bA=(A,\rhd,+,-)$. The category $\tr_{\bA}$ has as objects formal powers of $A$ and as arrows from
$A^m$ to $A^n$ relations $R:A^m\to A^n$ satisfying the following two conditions: 
if $\ba\, R\,\bb$ then for any element $c$ of $A$ 
\begin{itemize}
\item[TR(1)] $(c\rhd \ba)\, R\, (c\rhd \bb) $, and
\item[TR(2)] $\ba\rhd c =\bb\rhd c$.
\end{itemize}
It is an immediate consequence that the properties TR(1) and TR(2) of relations in $\tr_{\bA}$ are satisfied for any tuple $\bc$ in place of any element $c$: that is, if $R$ is a relation in $\tr_{\bA}$ and $\ba\, R\,\bb$ then for any tuple $\bc$ we have that $(\bc\rhd \ba)\, R\, (\bc\rhd \bb) $ and $\ba\rhd \bc =\bb\rhd \bc$.

Composition is the usual composition of relations. 

The monoidal structure of $\rel$ coming from the product of sets and of relations is easily seen to restrict to a strict monoidal structure on $\tr_{\bA}$.
\end{dfn} 

We next define a \emph{braiding} (\cite{JS93}) for the strict monoidal structure on $\tr_{\bA}$. This will require only the fact that tuples in an irack form a rack.
\begin{prop}
The functional relations
$$\tw_{m,n}:A^m\times A^n\to A^n\times A^m: (\ba,\bb)\mapsto (\ba\rhd\bb,\ba)$$ 
furnish $\tr_{\bA}$ with a braiding.
\end{prop}

\pf
We first check that $\tw_{m,n}$ is in $\tr_{\bA}$. The first property required is that for any  $c\in A$ we need that $\tw_{m,n}:c\rhd (\ba,\bb)\mapsto c\rhd (\ba\rhd\bb,\ba)$; that is, that
$$((c\rhd \ba)\rhd (c\rhd \bb),c\rhd\ba)=(c\rhd(\ba\rhd\bb),c\rhd\ba),$$ which follows by property R(3) of racks.
The second property required is that   $c\in A$ we need that $(\ba,\bb)\rhd c=  (\ba\rhd\bb,\ba)\rhd c$. That is, we need  that $\ba\rhd (\bb\rhd c)=  (\ba\rhd\bb)\rhd (\ba\rhd c)$, again property R(3) of racks.

The relation $(\ba,\bb)\mapsto (\bb,\ba\lhd\bb)$ is the inverse of $\tw_{m,n}$ by properties of racks.

Next to check that $\tw$ is natural. This amounts to two properties: if $R:A^m\to A^n$ is in $\tr_{\bA}$ then it is required that  $\tw_{n,p}(R\times 1_{A^p})=(1_{A^p}\times R)\tw_{m,p}$ and that $\tw_{p,n}( 1_{A^p}\times R)=(R\times 1_{A^p})\tw_{p,n}$. The first says that if $\ba\, R\, \bb$ and $\bc \in A^p$ then $\ba\rhd\bc=\bb\rhd\bc$. The second says that if $\ba\, R\, \bb$ and $\bc \in A^p$ then $\bc\rhd\ba \, R\, \bc\rhd\ba$.  These are the two properties TR(1),TR(2) required of $R$ to be in $\tr_{\bA}$.

Finally we need to check the two coherence conditions of twist. We will check just one, namely that
$$\tw_{m+n,p}=(\tw_{m,n}\times 1_{A^p})(1_{A^m}\times \tw_{n,p}):A^m\times A^n\times A^p\to A^p\times A^m\times A^n.$$
But under the first arrow $(\ba,\bb,\bc)$ goes to $((\ba,\bb)\rhd \bc,\ba,\bb)=(\ba\rhd(\bb\rhd \bc),\ba,\bb)$, and under the second arrow $(\ba,\bb,\bc)\mapsto (\ba,\bb\rhd\bc,\bb)\mapsto (\ba\rhd(\bb\rhd \bc),\ba,\bb)$.
\qed

We show now that the object $A$ of $\tr_{\bA}$ has a tangle algebra (\cite{FY89,RSW12}) structure.

\begin{prop}  The relations $\eta:I\to A\times A:()\sim (a,a^-)$ and $\varepsilon:A\times A\to I:(a,a^+)\sim ()$ furnish the object $A$ with a tangle algebra structure in $\tr_{\bA}$.
\end{prop}

\pf
To see that $\eta$ is in $\tr_{\bA}$ we need that $\eta:()\sim (c\rhd a,c\rhd a^-)$ which follows since
$c\rhd a^-=(c\rhd a)^-$.  We need also that $()\rhd c=(a,a^-)\rhd c$. Both sides evaluate as $c$. The proof that 
$\varepsilon$ is in $\tr_{\bA}$ is similar.

The fact the $\eta$ and $\varepsilon$ make $A$ self dual follows from the fact that $()^+$ and $()^-$ are inverse functions.

Finally the commutativity of $\eta$ and $\varepsilon$ both follow from the irack equation IR(2).
\qed

\section{Dirac's Belt Trick} 
In \cite{RSW12, MSW13} we proved that a blocked four torsion is equal to a blocked identity. We were not able to distinguish using $\tr_{\bG}$ for $\bG$ a group between a blocked two torsion and the identity and hence were unsure whether the blocked braid group on three strings has $6$ or $12$ elements. We are now able to present a proof of the unpublished result of Davide Maglia using the irack $\bA$ of example 1.2, namely that $\tr_{\bA}$ distinguishes between a blocked two torsion and a blocked identity on $3$ or more strings, and hence in particular that the blocked braid group on three strings has 12 elements.    

First let us consider the case of three strings. We will produce two relations $R:I\to A^3$ and $S:A^3\to I$ in $\tr_{\bA}$ such that if $T$ is the torsion on three strings then $SR$ is the empty set, while $ST^2R$ is the one point set.

The relation $R$ is $$R=\{(a,c,d),(b,e,f)\}$$ whereas the relation $S$ is $$S=\{(a,e,f),(b,c,d)\}.$$

A similar argument works for $n$ strings taking instead $$R=\{(a,c,d,1,1,\cdots ,1),(b,e,f,1,1,\cdots ,1)\}$$ and
$$S=\{(a,e,f,1,1,\cdots ,1),(b,c,d,1,1,\cdots ,1)\}.$$



\begin{thebibliography}{99}












\bibitem{CW59} John H. Conway, Gavin Wraith, unpublished correspondence, 1959.

%
%
\bibitem{RF92}   R.A. Fenn, C.P.Rourke,  \emph{Racks and Links in Codimension Two}, Journal of Knot Theory and its Ramifications 4, 343--406, 1992.









\bibitem{FY89} Peter Freyd and David Yetter, \emph{Braided compact monoidal categories with applications to low %
dimensional topology}, Adv. Math. 77, 156--182, 1989.






\bibitem {JS93} A. Joyal, R.H. Street, \emph{Braided monoidal
categories}, Adv. Math., Vol. 102, No. 1, 20--78, 1993.

\bibitem{J82} David Joyce, {\em A classifying invariant of knots: the knot quandle}, Journal of Pure and Applied Algebra 23: 37-65, 1982.



















\bibitem{MSW13} Davide Maglia, N. Sabadini, R.F.C. Walters, \emph{Blocked-braid Groups}, submmitted to Applied Categorical Structures, Springer Verlag, 2013.








 


\bibitem{RSW12} R. Rosebrugh, N. Sabadini, R.F.C. Walters, {\em Tangled Circuits}, Theory and Applications of Categories,  Vol. 26, No. 27, 743--767, 2012.






\bibitem{W} Gavin Wraith,  {\em A Personal Story about Knots}, http://www.wra1th.plus.com/gcw/rants/math/Rack.html, undated.


\end{thebibliography}
\end{document}